\documentclass[oneside,english]{amsart}
\usepackage[T1]{fontenc}
\usepackage[utf8]{luainputenc}
\usepackage{amsthm}
\usepackage{amssymb}
\usepackage{graphicx}
\usepackage[all]{xy}

\makeatletter
\numberwithin{equation}{section}
\numberwithin{figure}{section}
\theoremstyle{plain}
\newtheorem{thm}{\protect\theoremname}
  \theoremstyle{definition}
  \newtheorem{defn}[thm]{\protect\definitionname}
  \theoremstyle{remark}
  \newtheorem{rem}[thm]{\protect\remarkname}
  \theoremstyle{plain}
  \newtheorem{lem}[thm]{\protect\lemmaname}

\usepackage{amsfonts}\usepackage{amsthm}\usepackage{enumerate}\usepackage{mathrsfs}\usepackage{cancel}

\usepackage{MnSymbol}

\input xy\huge
\xyoption{all}

\DeclareFontFamily{OT1}{pzc}{}
\DeclareFontShape{OT1}{pzc}{m}{it}{<-> s * [1.10] pzcmi7t}{}
\DeclareMathAlphabet{\mathpzc}{OT1}{pzc}{m}{it}

\newcommand{\rr}{\mathbb{R}}

\usepackage{times}

\makeatother

\usepackage{babel}
  \providecommand{\definitionname}{Definition}
  \providecommand{\lemmaname}{Lemma}
  \providecommand{\remarkname}{Remark}
\providecommand{\theoremname}{Theorem}

\begin{document}

\title{Taylor Expansion Proof of the Matrix Tree Theorem - Part II}

\author{Amitai Zernik}

\address{Einstein Institute of Mathematics,}

\address{The Hebrew University of Jerusalem }

\address{Jerusalem, 91904, Israel}

\address{Email: amitai.zernik@gmail.com}
\begin{abstract}
The All Minors Matrix Tree Theorem states that the determinant of
any submatrix of a matrix whose columns sum to zero can be computed
as a sum over certain oriented forests. We offer a particularly short
proof of this result, which amounts to comparing Taylor series expansions. 
\end{abstract}
\maketitle

\section{\label{sec:Introduction}Introduction}

The following formulation of the All Minors Matrix Tree Theorem appears
in\textbf{ \cite{combinatorial-ammtt}.}
\begin{defn}
Let $U,W$ be subsets of $\{1,...,n\}$ of the same cardinality, $|U|=|W|=k$,
$k\geq1$. A \emph{forest from $U$ to $W$ }is an oriented graph
on $\{1,...,n\}$ which is a disjoint union of oriented trees with
the following properties:
\begin{enumerate}
\item Every tree in the forest contains exactly one vertex of $U$ and one
vertex of $W$.
\item The edges in each tree are oriented away from the vertex of the tree
belonging to $U$.
\end{enumerate}
We denote the set of all forests from $U$ to $W$ by $\mathbb{F}(U,W)$. 
\end{defn}
An $n\times n$ matrix $M$ is \emph{semi-laplacian} if all its columns
sum to zero. 
\begin{thm}
\label{thm:all-minors-mtt}Let $M$ be a semi-laplacian matrix. We
denote by $M(W,U)$ the submatrix of $M$ obtained by deleting the
$k$ rows indexed by $W$ and the $k$ columns indexed by $U$. Then
\begin{equation}
\det M(W,U)=\sum_{F\in\mathbb{F}(U,W)}\epsilon(U,W,F)A_{F}(M)\label{eq:All-Minors-MTT}
\end{equation}
where $A_{F}(M)=\prod_{(i,j)\in F}M_{ij}$ and $\epsilon(W,U,F)\in\{\pm1\}$
are signs given in Definition \ref{def:epsilon} below.\end{thm}
\begin{rem}
When $U=\{j\}$ and $W=\{i\}$ we obtain a direct generalization of
the well-known Matrix-Tree Theorem (which was discussed in the first
part of this paper). Namely, the determinant of the $i,j$ minor of
any semi-laplacian $M$ is expressed as a sum over all spanning trees
(since every spanning tree can be oriented in a unique way to obtain
a forest from $\{j\}$ to $\{i\}$).
\end{rem}
This theorem has been proven more than once. The proof given here
was guided by the same Taylor expansion yoga used in part 1 of the
paper. The signs are more subtle in this case, and were derived from
certain desirable cancellations which appeared in the proof (cf. Remark
\ref{rem:signs-from-proof}).

I would like to thank David Kazhdan, Nati Linial, Ori Parzanchevski,
Ron Peled, Ron Rosenthal and Ran Tessler for their suggestions and
comments. I am especially grateful to Nati Linial for suggesting reference
\textbf{\cite{combinatorial-ammtt}}.

\section{Proof of the All Minors Matrix Tree Theorem\label{sec:Proof-of-All-Minors-MTT}}

Throughout this section, $n$ will be some fixed positive integer.
On first reading, one may want to skip over Definition \ref{def:epsilon},
Lemma \ref{lem:epsilon-properties} and Remark \ref{rem:signs-from-proof}
and go directly to the proof of Theorem \ref{thm:all-minors-mtt},
so as not to lose sight of the forest for the signs.
\begin{defn}
\label{def:epsilon}Let $\mathcal{T}$ denote the set of all 3-tuples
$(U,W,F)$ where $U$ and $W$ are subsets of $\{1,...,n\}$ of the
same cardinality and $F$ is a forest from $U$ to $W$. Every such
$F$ defines a bijection $\pi_{F}:U\to W$. The sign of such a bijection,
$sgn(\pi)$, is defined to be the sign of the permutation that sorts
$(\pi(u_{1}),\pi(u_{2}),...,\pi(u_{k}))$ where $u_{1}<\cdots<u_{k}$
are the elements of $U$. We define $\epsilon:\mathcal{T}\to\{\pm1\}$
by
\[
\epsilon(U,W,F)=(-1)^{n+|U|}\cdot(-1)^{\sum_{u\in U}u+\sum_{w\in W}w}\cdot sgn(\pi_{F}).
\]

\end{defn}
The next lemma highlights the properties of $\epsilon$ that will
be used in the proof of Theorem \ref{thm:all-minors-mtt} below. 
\begin{lem}
\label{lem:epsilon-properties} $\epsilon:\mathcal{T}\to\{\pm1\}$
satisfies the following properties.

(a) $\epsilon(\{1,...,n\},\{1,...,n\},\emptyset)=+1$

(b) For any three elements $w_{0}\in W$, $i\not\in W$, $j\not\in U$
and any forest $F\in\mbox{Forests}(U\cup\{j\},W\cup\{i\})$ we have
\[
\epsilon'_{ij}(U,W)\epsilon(U\cup\{j\},W\cup\{i\},F)=\epsilon_{ij}''(F)\epsilon(U,W,F_{ij;w_{0}})
\]
where $\epsilon'_{ij}(U,W)=(-1)^{i+|\{w\in W|w<i\}|+j+|\{u\in U|u<j\}|}$
and $\epsilon''_{ij}(F)=-1$ if $i$ is a descendant of $j$ and $\epsilon''_{ij}(F)=+1$
otherwise.

(c) Let $i,j\in\{1,...,n\}$ and $w_{0}\in W$, and let $F$ be a
forest from $U$ to $W$ such that $F'=F\backslash\{(i,j)\}\cup\{(w_{0},j)\}$
is also a forest from $U$ to $W$. Then we have 
\[
\epsilon(U,W,F)=\epsilon(U,W,F').
\]

\end{lem}
We will prove Lemma \ref{lem:epsilon-properties} below. 
\begin{rem}
\label{rem:signs-from-proof}It is not hard to see that properties
(a) and (b) determine $\epsilon$ uniquely; indeed, we arrived at
Definition \ref{def:epsilon} by using property (c) to delete the
edges in $F$ one by one until the trivial forest is reached.\end{rem}
\begin{proof}
(of Theorem \ref{thm:all-minors-mtt}). When $|U|=|W|=n$ the claim
holds trivially since $\epsilon(U,W,\emptyset)=+1$ by part (a) of
Lemma \ref{lem:epsilon-properties}, so it is enough to prove the
theorem for $|U|=|W|=k$ given that it holds for $|U|=|W|=k+1$. 

Let $\mathcal{S}$ denote the vector space of semi-laplacian $n\times n$
matrices. Fix some $U,W\subset\{1,...,n\}$ with $|U|=|W|=k$, fix
$w_{0}\in W$, and let $\mathcal{D}:\mathcal{S}\to\rr$ and $\mathcal{F}:\mathcal{S}\to\rr$
denote the left and right hand sides of equation (\ref{eq:All-Minors-MTT}),
respectively. The tangent vectors $\frac{\partial}{\partial M_{ij}}-\frac{\partial}{\partial M_{w_{0}j}}$
form a basis for the tangent space to $\mathcal{S}$ at $M$. Since
$\mathcal{D}(0)=\mathcal{F}(0)=0$ and $\mathcal{S}$ is connected,
it is enough to prove 
\begin{equation}
\left(\frac{\partial}{\partial M_{ij}}-\frac{\partial}{\partial M_{w_{0}j}}\right)\mathcal{D}(M)=\left(\frac{\partial}{\partial M_{ij}}-\frac{\partial}{\partial M_{w_{0}j}}\right)\mathcal{F}(M)\label{eq:equate-ders}
\end{equation}
for all $i,j$ and $M\in\mathcal{S}$. Suppose first $i\not\in W$
and $j\not\in U$.

We have 
\begin{eqnarray}
(\frac{\partial}{\partial M_{ij}}-\frac{\partial}{\partial M_{w_{0}j}})\det M(W,U) & = & \frac{\partial}{\partial M_{ij}}\det M(W,U)\nonumber \\
 & = & \epsilon'_{ij}(U,W)\det M(W\cup\{i\},U\cup\{j\})\nonumber \\
 & = & \sum_{F\in\mathbb{F}(U\cup\{j\},W\cup\{i\})}\epsilon'_{ij}(U,W)\epsilon(U\cup\{j\},W\cup\{i\},F)A_{F}(M)\label{eq:minors-derivative}
\end{eqnarray}
where 
\[
\epsilon_{ij}'(U,W)=(-1)^{(i-|\{w\in W|w<i\}|)+(j-|\{u\in U|u<j\}|)}.
\]
and for the last equality we used our assumption that the claim holds
for $U,W$ of cardinality $k+1$.

On the right hand side of Eq \ref{eq:equate-ders}, we have 
\[
(\frac{\partial}{\partial M_{ij}}-\frac{\partial}{\partial M_{w_{0}j}})\sum_{F}\epsilon(U,W,F)A_{F}(M)
\]

It is not hard to see that there are pairs of forests $F,F'$ with
$\frac{\partial}{\partial M_{ij}}A_{F}(M)=\frac{\partial}{\partial M_{w_{0}j}}A_{F'}(M)$;
as we shall see, these forests cancel each other's contribution, and
the remaining contributions can be interpreted as coming from forests
from $U\cup\{j\}$ to $W\cup\{i\}$. We now make this precise.

We have

\[
\frac{\partial}{\partial M_{ij}}A_{F}(M)=\begin{cases}
A_{F\backslash\{(i,j)\}}(M) & \mbox{if }(i,j)\in F\\
0 & \mbox{otherwise}
\end{cases}
\]
where $F\backslash\{(i,j)\}$ is the forest obtained by deleting the
edge from $i$ to $j$, similarly for $\frac{\partial}{\partial M_{w_{0}j}}A_{F}(M)$.
Suppose $(i,j)\in F$. $F\backslash\{(i,j)\}$ can also be written
as $F'\backslash\{(w_{0},j)\}$ for the oriented \emph{graph} $F'=F\backslash\{(i,j)\}\cup\{(w_{0},j)\}$,
in which the vertex $j$ has been ``reattached'' to $w_{0}$. $F'$
will be a forest from $U$ to $W$ if and only if there's no oriented
path from $j$ to $W$ in $F$. 

\begin{figure}
\includegraphics[scale=0.3]{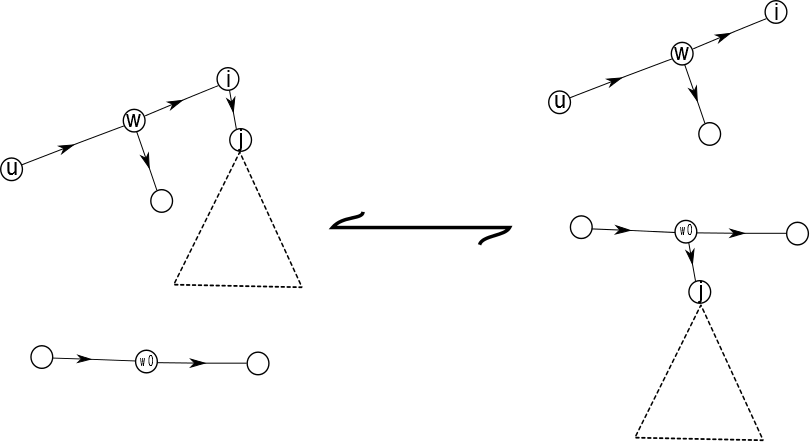}

\caption{\label{fig:Two-cancelling-forests}Two cancelling forests, $F$ and
$F'$. The dashed triangle represents the subtree consisting of $j$'s
descendants, which is disconnected from $i$ and reattached to $w_{0}$.}
\end{figure}

Conversely, we can start with any $F'$ such that $(w_{0},j)\in F'$
and consider the oriented graph $F=F'\backslash\{(w_{0},j)\}\cup\{(i,j)\}$.
$F$ will be a forest from $U$ to $W$ if and only if there's no
oriented path from $ $$j$ to $i$ in $F'$. See Fig. \ref{fig:Two-cancelling-forests}.
The operations $F\mapsto F'$ and $F'\mapsto F$ define a bijection
between the two subsets of $\mathbb{F}(U,W)$:
\[
\left\{ \mbox{\ensuremath{(i,j)\in F}and \ensuremath{\not\exists}path }j\to W\mbox{ in \ensuremath{F}}\right\} \simeq\left\{ \mbox{\ensuremath{(w_{0},j)\in F'}and \ensuremath{\not\exists}path \ensuremath{j\to i}in \ensuremath{F}'}\right\} .
\]
By part (c) of Lemma \ref{lem:epsilon-properties} $\epsilon(U,W,F)=\epsilon(U,W,F')$,
so these terms cancel in pairs in $(\frac{\partial}{\partial M_{ij}}-\frac{\partial}{\partial M_{w_{0}j}})\sum_{F}\epsilon(U,W,F)A_{F}(M)$.
The other contributing forests contain either $(i,j)$ or $(w_{0},j)$,
but are such that if we delete the incoming edge to $j$ we obtain
a forest from $U\cup\{j\}$ to $W\cup\{i\}$, see Fig. \ref{fig:non-cancelling}.
Conversely, given a forest $F$ from $U\cup\{j\}$ to $W\cup\{i\}$
we can add precisely one of the edges $(w_{0},j)$ or $(i,j)$ to
obtain a forest, denoted $F_{ij;w_{0}},$ from $U$ to $W$. To see
this, consider the unique oriented path from $j$ to $W$; if it hits
$i$, we must take $F_{ij;w_{0}}=F\cup\{(w_{0},j)\}$, otherwise $F_{ij;w_{0}}=F\cup\{(i,j)\}$.
We call the operation $F\mapsto F_{ij;w_{0}}$ \emph{a gluing} (of
$j$).

\begin{figure}
\includegraphics[scale=0.3]{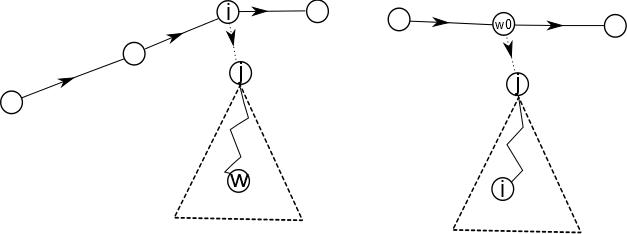}

\caption{\label{fig:non-cancelling}Two types of forests (only one tree is
shown of each) which do not cancel with any other forest. In both
cases, if we delete $j$'s incoming edge we obtain a forest from $U\cup\{j\}$
to $W\cup\{i\}$.}
\end{figure}

The above discussion shows that

\begin{equation}
(\frac{\partial}{\partial M_{ij}}-\frac{\partial}{\partial M_{w_{0}j}})\sum_{F}\epsilon(U,W,F)A_{F}(M)=\sum_{F\in\mathbb{F}(U\cup\{j\},W\cup\{i\})}\epsilon_{ij}''(F)\epsilon(U,W,F_{ij;w_{0}})A_{F}(M)\label{eq:forests-derivative-take2}
\end{equation}
where $\epsilon''_{ij}(F)=-1$ if $i$ is a descendant of $j$ and
$\epsilon''_{ij}(F)=+1$ otherwise. Comparing eq (\ref{eq:forests-derivative-take2})
and eq (\ref{eq:minors-derivative}) and using property (b) of Lemma
\ref{lem:epsilon-properties}, eq (\ref{eq:equate-ders}) is proven
for $j\not\in U$ and $i\not\in W$. 

Finally, observe that if $j\in U$ or $i\in W$ then $\left(\frac{\partial}{\partial M_{ij}}-\frac{\partial}{\partial M_{w_{0}j}}\right)\mathcal{D}(M)=0$,
and so we need to show that the RHS of eq (\ref{eq:equate-ders})
vanishes. If $j\in U$ then $\frac{\partial}{\partial M_{ij}}\mathcal{\mathcal{F}}(M)=\frac{\partial}{\partial M_{w_{0}j}}\mathcal{\mathcal{F}}(M)=0$
since $U$ has no incoming edges. If $i\in W$ then \emph{all} the
terms of $(\frac{\partial}{\partial M_{ij}}-\frac{\partial}{\partial M_{w_{0}j}})\sum_{F}\epsilon(U,W,F)A_{F}(M)$
cancel in pairs. E.g., if $(i,j)\in F$ then there cannot be a path
from $j$ to $W$, because that would imply there are two elements
of $W$ in the same tree, and similarly for $(w_{0},j)\in F$.
\end{proof}

\begin{proof}
(of Lemma \ref{lem:epsilon-properties}). 

Property (a) is immediate. 

For property (b), assume first that $\epsilon''(F)=+1$. That is,
$i$ is not a descendant of $j$. Consider the following diagram:

\[
\xymatrix{\ar@/{}_{1pc}/@{.>}[rrrrr] & *+[o][F]{u_{1}} & *+[o][F]{u_{2}} &  & \cdots & *+[o][F]{j}\ar[dll] & \cdots & *+[o][F]{u_{r}}\ar[dl]\ar@{-->}[dllll] & \cdots\\
\ar@{.>}[u] & *+[o][F]{w_{1}} & \cdots & *+[o][F]{w_{s}} & \cdots &  & *+[o][F]{i}\ar@/_{1pc}/@{.>}[llllll] & \cdots
}
\]

The full arrows represent part of the bijection $\pi_{F}$ for the
forest $F$; we assume the $u$'s and $w$'s are sorted on the top
and bottom row, respectively; in such a diagram the sign of the bijection
is given by the parity of the number of arrow intersections, and is
independent of how the arrows are drawn%
\footnote{Strictly speaking, the arrows must be drawn inside the strip and must
intersect transversally for this to hold.%
}. The bijection $\pi_{F_{ij;w_{0}}}$ for the glued forest is obtained
by erasing the full arrows, together with $j$ and $i$, and replacing
them with the dashed arrow. 

Since the dashed arrow can be drawn by tracing the full arrow from
$u_{r}$ to $i$, the dotted path, and then the full arrow from $j$
to $w_{s}$, we see that the sign difference $sgn(\pi_{F})sgn(\pi_{F_{ij;w_{0}}})$
is given by the parity of the number of intersection of the dotted
path with the other arrows in the diagram (i.e., those \emph{not}
shown) plus the number of intersections of the full arrows with each
other. In the case drawn, where $j<u_{r}$ and $w_{s}<i$, the full
arrows do not intersect; the top dotted arrow intersects the outoing
arrows from the set of vertices $\{u\in U:u<j\}$; the bottom dotted
arrow intersects the incoming arrows to $\{w\in W:w<i\}$, except
$w_{s}$'s incoming arrow is a full arrow, so we subtract one from
the count, to find that 
\[
sgn(\pi_{F})sgn(\pi_{F_{ij;w_{0}}})=(-1)^{|\{u\in U:u<j\}|+|\{w\in W:w<i\}|-1}
\]
in this case. It is not hard to see that the same formula holds for
the three other cases: $j<u_{r}$ and $i<w_{s}$; $u_{r}<j$ and $w_{s}<i$;
$u_{r}<j$ and $i<w_{s}$. The verification of property (b) for $\epsilon''(F)=+1$
now follows by a straightforward computation.

The case $\epsilon''(F)=-1$, when $i$ is a descendant of $j$, is
similar but simpler, since there's a single arrow from $i$ to $j$
that needs to be discarded. The sign difference in this case is

\[
sgn(\pi_{F})sgn(\pi_{F_{ij;w_{0}}})=(-1)^{\{u\in U:u<j\}+\{w\in W:w<i\}}
\]
and the result follows. The proof of property (b) is complete.

For property (c), observe that $F,F'$ are both forests from $U$
to $W$ with the same set of oriented paths from $U$ to $W$, so
$\pi:U\to W$ is also the same.
\end{proof}
\bibliographystyle{plain}
\bibliography{matrix_tree_bibilography}

\end{document}